\documentstyle{article}

\title{Breather Solutions of the Nonlinear Wave Equation}
\author{Satyanad Kichenassamy\thanks{Supported by the N.S.F.
            under contract \# DMS-8504033 and by the Army.}\ \thanks{To appear in Communications on Pure and Applied Mathematics.}\\
       Courant Institute of Mathematical Sciences}
\date{{\bf Appeared in}: Communications in Pure and Applied Mathematics
{\bf 44} (1991) 789--818.}

\newtheorem{theorem}{Theorem}

\begin{document}

\maketitle
\tableofcontents

\section{Introduction.}

  \subsection{Generalities.}
    This paper is devoted to the problem:
  \begin{equation}
    \left\{     \begin{array}{l}
      u_{tt} - u_{xx} + g(u) = 0  \\
      u(x,t+T) = u(x,t) ; \;
      u(x,t) \rightarrow 0 \mbox{ as } |x| \rightarrow \infty ; \;
      u_{t} \not\equiv 0.
        \end{array} \right.
  \end{equation}

    A solution of this problem is called a {\em breather}.

    Throughout this paper, $g$ is assumed to be entire with $g(0)=0$
and $g'(0) = 1$. We write
\[ g(u) = \sum_{m=1}^{\infty} g_{m} u^{m}. \]

 When
$g(u) = \sin u$, one can write down an explicit breather for every period $T$
greater than $2 \pi$; it reads:
  \begin{equation}      \label{SGB}
    u_{SG} = 4 \arctan (\frac{\varepsilon}{\sqrt{1-\varepsilon^{2}}}
      \frac{\cos(t \sqrt{1-\varepsilon^{2}})}
           {\cosh(\varepsilon x)}),
  \end{equation}
where $\varepsilon$ is related to $T$ by \[ \varepsilon^{2}+
(2\pi /T)^{2} = 1. \] This is connected to the Inverse Scattering Transform.
Some solitons are therefore breathers.

    Observe that (\ref{SGB})
is given  as a series in $\varepsilon$,
$e^{-\varepsilon x}$ and $\cos (t\sqrt{1-\varepsilon^{2}})$,
convergent on any domain of the form
\[ \{|\varepsilon| \leq \frac{1}{\sqrt{2}} + \delta ;
     \varepsilon x \geq A(\delta) \} \]
with a suitable function $A(\delta)$.

    We sketch in an appendix some previous attempts at proving or
disproving the existence of breathers, and give selected references.

    We now list the results of the present work. Various improvements
and other approaches will be considered in subsequent papers.

  \subsection{Results.}

    Our first two results
(Theorems 1 and 2) show that one can construct and characterize two different
types of {\em formal solutions} to our problem.

     One is of the type
  \begin{equation}      \label{series1}
         \sum_{k=1}^{\infty} \left( \varepsilon^{k}
         \sum_{\begin{array}{c}
        0 \leq q \leq l \leq k ; \,\, l \geq 1\\
        k \equiv q \equiv l \pmod{2}
          \end{array}}
      a_{kql} (1/\cosh(\varepsilon x))^{l}
                  \cos(qt\sqrt{1-\varepsilon^{2}}) \right)
  \end{equation}
 in powers of the parameter $\varepsilon$ defined above.

    The other has the form
  \begin{equation}          \label{series2}
    \sum_{l=1}^{\infty}  u_{l} (t,\varepsilon)e^{-l\varepsilon x},
  \end{equation}
in powers of
$e^{-\varepsilon x}$. These series are constructed in Sections 2 and 3.

    We shall prove that the
coefficients of the second series are polynomials in
$\cos (t\sqrt{1-\varepsilon^{2}})$,
 but are {\em rational} in $\varepsilon$.

    Theorem 3, stated and proved in Section 4, shows that (\ref{series1})
 {\em does not} define an analytic function of $\varepsilon$,
$e^{-\varepsilon x}$ and $\cos (t\sqrt{1-\varepsilon^{2}})$
 on any domain of the form
\[ \{ |\varepsilon| \leq \frac{1}{\sqrt{2}} + \delta ; x \geq A \} \]
{\em unless} $g(u)$ has one of the three forms:
\[  C \sin (\alpha u), \]
\[  C \sinh (\alpha u), \]
\[  Cu.  \]
Only the first gives rise to breathers, which are merely scaled versions of
(\ref{SGB}). The sine-Gordon equation is thereby singled out among
non-linear wave equations.



    The proof of Theorem 3 involves a detailed study of the series
(\ref{series2}). We have in that connection written a program in
``Mathematica'' which generates this series. The program was run on a
Sun 3 of the Academic Computing Facility at New York University.
The results of this calculation
reveal that the terms of this series, as functions of the coefficients
of $g$, are very strongly related to one another. We give in Section 4
some first results in that direction.

    Theorem 4, stated and proved in Section 5, shows that the series
  (\ref{series2})      {\em does converge} for every fixed
{\em real} value of $\varepsilon$, and large $x$,
{\em provided} that $g(u) = \sum_{m \;\mbox{\scriptsize odd}} g_{m} u^{m}$ with
\[   |g_{m}| \leq \alpha ^{(m-1)} / m!\]
for every odd $m$ and some positive $\alpha$.
More prcisely, convergence takes place when $\varepsilon x$ lies in the domain
of convergence of the corresponding series for
$g(u) = \sinh(\alpha u)/\alpha$, which reads
\[ u_{SHG} :=  \frac{4i}{\alpha}
   \arctan \left( \frac{\varepsilon}{\sqrt{1-\varepsilon^{2}}}
                  \frac{\cos (t\sqrt{1-\varepsilon^{2}})}
                       {i \sinh(\varepsilon x)}
           \right).                    \label{SHG}                          \]

    This does not contradict the former result since convergence only
takes place for {\em large} $x$.

    When $g$ is not odd, we cannot expect a similar result,
since the functions $u_{l}(t,\varepsilon)$ have in general {\em real}
poles for $\varepsilon^{2} = 1/m^{2}$, $m$ an integer --- as opposed to
pure imaginary poles for odd $g$.

    We can draw the following consequences from these results:
  \begin{enumerate}
    \item The sine function does have a special status in this
problem, and for reasons not directly related to the
Inverse Scattering Transform or B\"{a}cklund transformations.
The introduction of $\varepsilon$ as a new variable creates a degeneracy
which prevents the existence of a solution holomorphic in the
three variables $\varepsilon$, $t$, and $e^{-\varepsilon x}$.
    \item Our convergence result means that when $\varepsilon$ is kept
fixed, there is in some cases of interest a {\em convergent} series solution
for $x$ large enough, while it would have seemed difficult, in view of the
divergence result, to expect more than a series of asymptotic type.
  \end{enumerate}

    Theorem 5, stated and proved in Section 6, shows that for any
odd $g$, without growth conditions on the coefficients, there always
exists a solution of our equation which tends to zero as $x$ tends
to $+\infty$ only, and which belongs to every $H^{s}$ class in
$x$ and $t$. This follows from a somewhat detailed study of a
stable manifold-type argument. Our treatment differs from
others mentioned in the Appendix in that we establish period-independent
estimates and actually single out one solution in a
multi-dimensional stable manifold.

\vspace{1 ex}

    {\bf Acknowledgements.} I should like to express my thanks to
Peter D. Lax, for suggesting this problem, and for his constant
encouragement and his very helpful remarks. I am also glad to thank
L. Nirenberg for an important observation on a first draft,
and H. Br\'{e}zis, P.A. Deift and H.P. McKean for fruitful discussions.

\section{First Formal Solution.}

\subsection{Statement of the result.}

    We shall construct in the present section a formal solution to
the equation:
\begin{equation}
  u_{tt}-u_{xx}+ u + \sum_{m=2}^{\infty} g_{m} u^{m} = 0.
\end{equation}

    We seek the solution $u$ in the form
\begin{equation}
  u(x,t) = v(\varepsilon x, \sqrt{1-\varepsilon^{2}} \, t, \varepsilon)
\end{equation}
where $\varepsilon$ is a new parameter, and $v$ is a formal series:
\begin{equation}
    v = \sum_{k=1}^{\infty} \varepsilon^{k} v_{k}.
\end{equation}

    Letting $\xi:= \varepsilon x$ and
$\tau := \sqrt{1-\varepsilon^{2}} \, t$, we find that $v$ satisfies the
following equation:
\begin{equation}        \label{scaledeq}
v_{\tau\tau} + v -\varepsilon^{2}(v_{\tau\tau}+v_{\xi\xi})
+\sum_{m=2}^{\infty} g_{m} v^{m} = 0.
\end{equation}



    Our result is that under a simple condition on $g$,
there are precisely two nonzero such series,
if we require the functions $v_{k}$ to be even in their arguments
$\xi$ and $\tau $,
decaying in $\xi$, and of period $2\pi/\sqrt{1-\varepsilon^{2}}$ in $t$:
  \begin{theorem}
    Assume that
\begin{equation}
  \lambda := \frac{5}{6} g_{2}^{2} - \frac{3}{4} g_{3} > 0.
\end{equation}
Then (\ref{scaledeq}) has a formal solution of the form
\[ \sum_{k=1}^{\infty} \varepsilon^{k} v_{k}(\tau,\xi).  \]
There are only two nonzero such series, determined by the
following conditions:
\begin{enumerate}
    \item  all $v_{k}$ are $2\pi$ periodic in $\tau$ and tend to
zero as $\xi \rightarrow \pm \infty$;
    \item  all $v_{k}$ are even in $\tau$ and $\xi$.
\end{enumerate}
One then finds that
\begin{enumerate}
\item
\(   v_{1} = \pm S \cos \tau \mbox{ with }
              S = \frac{\sqrt{2/\lambda}}{\cosh \xi}. \)
\item The functions $v_{k}$ have the form
\[ v_{k} = \sum_{0 \leq q\leq k} a_{kq}(S) \cos q\tau  \]
where $a_{kq}$ is a polynomial of degree $\leq k$, vanishes
if $k+q$ is odd, has the same parity as $k$, and is always divisible by $S$.
\end{enumerate}
  \end{theorem}

{\em Remarks:}
\begin{enumerate}
    \item We have thus obtained and characterized the series (\ref{series1})
 announced in the Introduction.
    \item $\lambda > 0$ for $g = u - u^{3}$ and for
$g = (-(1+u)+(1+u)^{3})/2$ (``$\varphi^{4}$-model;'' the factor $1/2$
makes $g'(0)$ equal to one.).
\end{enumerate}

\subsection{Proof of Theorem 1.}

    The strategy is a follows:

    We first compute the first three terms of our series, which will
contain two arbitrary functions (\S 2.2.1). The computation of $v_{3}$
will impose the value of $v_{1}$. We next show (\S 2.2.2)
 that the linearization
of the equation determining $v_{1}$ has an inverse, mapping polynomials
of the form $S^{3} P(S^{2})$ to polynomials of the form $S Q(S^{2})$.
This will enable us to construct the coefficients $v_{k}$
inductively (\S 2.2.3).

  \subsubsection{First three terms of the series.}

 Substitution of the series in (\ref{scaledeq}) results in the following:
\[  v_{1\tau\tau} + v_{1} = 0;      \]
\[  v_{2\tau\tau} + v_{2} + g_{2} v_{1}^{2} = 0;    \]
\[  v_{3\tau\tau} + v_{3} + 2g_{2}v_{1}v_{2} + g_{3}v_{1}^{3}-\Delta
                       v_{1}    = 0,            \]
where $\Delta := \partial_{\tau\tau} + \partial_{\xi\xi}$.

    Since the functions $v_{k}$ are even in both $\tau$ and $\xi$,
it follows that there are functions $S(\xi)$,$\sigma_{2}(\xi)$,
even and decaying as $|\xi| \rightarrow \infty$, such that:
\[  v_{1} = S(\xi) \cos \tau;           \]
\[  v_{2} = -(1/2) g_{2} S^{2} (1-(1/3) \cos 2\tau)
        +\sigma_{2}\, \cos \tau;            \]
\[  v_{3\tau\tau} + v_{3} =
(S''-S+\lambda S^{3}) \cos \tau -(g_{2}^{2}/6+g_{3}/4)S^{3}\cos3\tau
          - 2 g_2 \sigma_2 S \cos^2 \tau.\]
Here, $\lambda$ is the quantity defined in Theorem 1 and the prime
stands for differentiation with respect to $\xi$.

    Periodicity of $v_{3}$ requires that $v_{3\tau\tau} + v_{3}$
be orthogonal to $\cos \tau$, which imposes
\[  S'' - S +\lambda S^{3} = 0.         \]

    This equation has, up to  translation, exactly two nontrivial solutions tending to zero
as $|\xi| \rightarrow \infty$, namely:
\[  \pm \frac{\sqrt{2/\lambda}}{\cosh \xi}.     \]

    They are clearly even. Each of these will generate a formal solution to
our problem.

    We fix $S$ to be one of these two functions from now on.

    $v_{3}$ is then given by:
\[  v_{3} = \frac{S^{3}}{8}(g_{2}^{2}/6 + g_{3}/4) \cos 3\tau
          - g_2 \sigma_2 S (1-\frac{1}{3}\cos 2\tau) + \sigma_{3} \cos \tau,        \]
where $\sigma_{3}$, like $\sigma_{2}$ is still unknown. Here again, as $v_{3}$
is even in $\tau$, no term in $\sin\tau$ appears.

  \subsubsection{An auxiliary Schr\"{o}dinger operator.}

    We here give some properties of the linearization of
\[  S'' - S + \lambda S^{3} = 0 \]
at $S_{0} = \pm\frac{\sqrt{2/\lambda}}{\cosh \xi}$. This linearization is the
Schr\"{o}dinger operator
\[  L := \partial_{\xi}^{2} - 1 +\frac{6}{\cosh^{2} x}, \]
which is a translate of the  Schr\"{o}dinger equation with ``2-soliton''
 potential (one of the ``Bargmann potentials''). As an operator on
$L^{2}$, it possesses two eigenvalues, one of them being zero.
The other is $3$, with eigenfunction $S^{2}$.
$dS/d\xi$ is clearly an eigenfunction for the eigenvalue zero, and it is
simple, since decaying potentials are in the limit-point
case. The essential spectrum is the half-line $(-\infty,-1]$. It
easily follows that $L$, on the space of even, square-summable
functions is {\em invertible}.

    For the proof of Theorem 1, we shall need a more precise result:

\vspace{2 ex}

{\bf Lemma:} {\em Let P be a polynomial. The only even, square-summable
solution of $L\sigma = S^{3} P(S^{2})$ has the form
$\sigma = S Q(S^{2})$, where $Q$ is a polynomial.
It is also the only even solution tending to zero
as $\xi$ tends to infinity. }

{\em Proof:} It suffices to find one solution of the required form, since we
know that $L$ is invertible on even functions. Now, for any integer $p\geq0$,
one has
\[     L(S^{p}) = S^{p}(p^{2} - 1 + (3 - p(p+1)/2)\lambda S^{2})    \]
which implies the result.



\vspace{2 ex}

{\em Remarks:}
\begin{enumerate}
    \item  There is a very simple way to use the preceding calculation
to prove the invertibility of $L$ on even functions. As the essential
spectrum is bounded away from zero, it suffices to show that 0 is not an
eigenvalue for $L$ on even, square-summable functions. But if there
were a nonzero $\sigma$ such that $L\sigma = 0$, $\sigma '(0) = 0$, then

\begin{enumerate}
  \item $\sigma$
has a finite number of zeros (because $\sigma''$ and $\sigma$ have the same sign
far out), and none of them is 0;
  \item if $\alpha_{1}, \ldots, \alpha_{r}$ are the positive zeros of $\sigma$,
in increasing order, and if, say, $\sigma(0) > 0$, then
\[   \Sigma :=  S^{3}(S^{2} -S^{2}(\alpha_{1})) \ldots
              (S^{2} - S^{2}(\alpha_{r})) \]
has the same sign as $\sigma$, since $S$ is decreasing; it lies in the range
of $L$, by the Lemma.
\end{enumerate}
    But $\sigma$ and $\Sigma$ should then be
$L^{2}$-orthogonal, which is impossible. This shows that such a $\sigma$ cannot
exist.
    \item  An analogous argument shows that
\[  \sigma'' - \sigma + \beta S^{2} \sigma = 0  \]
has no even, square-integrable solution if $\beta \neq p(p+1)\lambda/2$
for every odd integer $p$. (Compare \cite[p. 60]{ince}.)
\end{enumerate}

  \subsubsection{Induction and end of proof.}

    The argument is in six Steps.

\vspace{1 ex}

    {\em Step 1: Induction hypothesis.}

    Assume, by induction, that we have, for some $m \geq 4$, found
$v_{1} \ldots v_{m-3}$, and that they have the form:
\begin{equation}
v_{k} = \sum_{\begin{array}{c}
        0 \leq q \leq k \\
                 q \equiv k \pmod{2}
              \end{array}}
a_{kq}(S) \cos q\tau,
\end{equation}
where the $a_{kq}$ are polynomials of degree $k$ at most, of the same
parity as $k$, and divisible by $S$: assume also that the relations
obtained by setting to zero the coefficients of
$\varepsilon, \ldots, \varepsilon^{m-1}$ in (\ref{scaledeq})
determine $v_{m-2}$ and $v_{m-1}$ as well, in the form
\[ \sum_{q \neq 1} a_{kq}(S) \cos q\tau + \sigma_{k} \cos \tau \]
where $k=m-1$ or $m-2$, the $a_{kq}$ are polynomials with the properties
listed above, and $\sigma_{m-2}$, $\sigma_{m-1}$ are arbitrary
functions of $\xi$.

    The remainder of the present Section is devoted to showing that
the above properties hold for $m=4$, and that if they hold for some $m$,
they also hold for $m+1$. This will prove Theorem 1.

\vspace{1 ex}

    {\em Step 2: Hypothesis holds for $m=4$.}

This follows from the expressions for $v_{1}$, $v_{2}$ and $v_{3}$ found
above.

\vspace{1 ex}

    {\em Step 3: Equation for $v_{m}$.}

    Setting to zero the coefficient of $v_{m}$ in (\ref{scaledeq}),
we find that $v_{m}$ solves
\begin{eqnarray}        \label{vmeqn}
  v_{m\tau\tau} + v_{m} & = & \Delta v_{m-2}
                -2g_{2}(v_{1}v_{m-1}+v_{2}v_{m-2})
             \nonumber \\
                        &   & -3g_{3}(v_{1}^{2}v_{m-2})+\phi_{m}
                                      (v_{1},\ldots,v_{m-3}).
\end{eqnarray}

    In this equation, $\phi_{m}$ is a polynomial in its arguments.
>From the induction hypothesis, the right-hand side of (\ref{vmeqn}) is a
trigonometric polynomial, and we must ascertain

    i) that the coefficient of $\cos j\tau$ is, for every $j \neq 1$,
a polynomial in $S$ of the right form,

    ii) that the coefficient of $\cos \tau$ is zero, so that no
secular term arises.

    It is the second condition that determines $\sigma_{m-2}$.
As all we know from (\ref{vmeqn}) is the value of $v_{m\tau\tau} + v_{m}$,
$v_{m}$ will be determined up to the addition of a term
$\sigma_{m}(\xi) \cos \tau$, which will be computed only in the calculation of
$v_{m+2}$.

    We examine the above two points in order, thereby finishing the proof
of the Theorem.

\vspace{1 ex}

    {\em Step 4: Examination of the coefficient of
$\cos j\tau$ for $j \neq 1$.}

    We show here that the right-hand side of (\ref{vmeqn}) contains
only powers of $S$ of the same parity as $m$ and that it contains
$S^{2}$ (resp. $S$) as a factor for $m$ even (resp. odd):

    All terms in the right-hand side of (\ref{vmeqn}) come either
from $\Delta v_{m-2}$ or from the expansion of some product of the form
\begin{equation}    \label{s}
v_{r_{1}}^{\alpha_{1}} \ldots v_{r_{s}}^{\alpha_{s}}
\end{equation}
where
\[ s\geq 1,  \; \alpha_{i} \geq 1, \;
    \sum_{i} \alpha_{i}r_{i} = m, \; 2 \leq \sum \alpha_{i} \leq m. \]



    Under the induction hypothesis, such a
product can only be a combination of terms of the form $S^{q}\cos l\tau$,
where $q$ and $l$ have the same parity as $\sum_{i} \alpha_{i} r_{i}$,
that is, as $m$ itself.

    On the other hand, the $\Delta v_{m-2}$ term has the required
form, since
\[ (S^{m})'' = S^{m}(m^{2} - m(m+1)\lambda S^{2}/2).    \]

    Finally, we must check the powers of $S$ which occur:

    $\Delta v_{m-2}$  can contain $S^{q} \cos l\tau$ only if $q \geq 1$
(in case $l$ is odd) or $q\geq 2$ (in case $l$ is even)
because of the induction hypothesis.

    For the remaining terms, we shall  distinguish three cases
according to the number $s$ of factors in (\ref{s}):
\begin{enumerate}
    \item $s=1$: This corresponds to the term $-g_{m}v_{1}^{m}$, which has
a factor $S^{m}$. As $m \geq 2$, this term has the required form.
    \item $s=2$: This corresponds to terms $v_{k}v_{m-k}$. If $m$ is
even, they have $S^{2}$ as a factor, by the induction hypothesis. If $m$ is odd,
either $k$ or $m-k$ is odd, and the other even; these terms therefore
contain at least $S.S^{2}$ as a factor.
    \item $s \geq 3$: A term generated by such a product contains at
least $S$ to the power $\alpha_{1} + \cdots + \alpha_{s}$ which is
greater than or equal to
$s$, and therefore greater than or equal to $3$: This term contains $S^{3}$ as a factor.
\end{enumerate}

    {\em Step 5: Examination of the coefficient of $\cos \tau$.}

    The equation expressing the vanishing of the coefficient of
$\cos \tau$ in the right-hand side of (\ref{vmeqn}) reads
\begin{equation}
 L\sigma_{m-2} = S^{3}\varphi_{m}(S^{2})
\end{equation}
where $L$ is the operator studied in Section 2.2, and $\varphi_{m}$
is a polynomial of degree less than or equal to $ (m-3)/2$ obtained by collecting terms in
$\cos \tau$ in the expansion of $\phi_{m}$. $L\sigma_{m-2}$ comes from
the substitution of $v_{1}$ and $v_{2}$ in (\ref{vmeqn}). Note that
$v_{m\tau\tau} + v_{m}$, as given by (\ref{vmeqn}), contains only even
cosines for even $m$; we must therefore have $\varphi_{m} = 0$ for even $m$.

    It now follows  from the results of \S 2.2 that (13)
has a unique even and decaying solution, which has the
required form.

 \vspace{1 ex}

    {\em Step 6: End of proof.}

    We have seen that the induction defined above generates the
coefficients $v_{k}$ recursively, and that they are uniquely
determined by the assumptions of decay and parity. The proof is now
complete.

  \subsection{Other means of generating the formal series.}

    We briefly mention here two other means of constructing a formal
solution in increasing powers of $\varepsilon$:

\vspace{1 ex}

    (a)  Solve the Cauchy problem for (\ref{scaledeq}) with data
$(S(\xi),0)$, $S$ even, decaying. Let
\[      T(S,\varepsilon) = v_{\tau}(\pi)/\varepsilon^{2}.   \]
Then a direct calculation leads to
\[  T(S,0) = \mbox{const.} (S'' - S +\lambda S^{3}),    \]
so that $T_{S}(S,0)$ is invertible, by the results of \S 2.2. The existence of
a formal solution to $T(S,\varepsilon) = 0$ in powers of $\varepsilon$
follows easily.

\vspace{1 ex}

    (b)  Let $F(v,\varepsilon)$ be the right-hand side of (\ref{scaledeq}).
Let $Q$ denote the projection on functions orthogonal to $\cos \tau$, $P$
being the complementary projection. Write
\[  v(\xi,\tau) = S(\xi) \cos \tau + \varepsilon w(\xi,\tau),   \]
where
\[   \int_{0}^{2\pi} w \, \cos \tau \, d\tau = 0.       \]
One can then see that the operator
\[  \tilde{F}(v,\varepsilon) = QF/\varepsilon^{2} + PF/\varepsilon  \]
is well-defined.  We consider it as a function of the variables
$S$, $w$, and $\varepsilon$. When $\varepsilon$ equals $0$, $\tilde{F} = 0$
reduces to
\[   \left\{
\begin{array}{l}
S'' - S +\lambda S^{3} = 0  \\
w_{\tau\tau} + w = - g_{2} S^{2} \cos^{2} \tau,
\end{array}
\right.     \]
as it should.

    Once again the linearization of $\tilde{F} = 0$ when $\varepsilon=0$
is invertible at the nonzero solutions, and the existence of a formal solution
in powers of $\varepsilon$ follows. This second approach was suggested
by techniques used for the water wave problem (see \cite{beale}).

    Both approaches lead to a formal series by the same process, which may
be worth mentioning.    When an equation takes the form
$T(u,\varepsilon) =0$, where $T(u_{0},0) = 0$ and $T_{u}(u_{0},0)$ is
invertible in some sense, then, under fairly general conditions, one can
show that $T(u,\varepsilon)=0$ has a formal solution in increasing powers of
$\varepsilon$. One can put this observation in a more precise form, which
however will not be needed in the present paper.

\section{Second Formal Solution.}

  \subsection{Statement of the result.}

    The preceding section has shown that one can achieve a formal
solution of our problem in increasing powers of $\varepsilon$, and that
the series thus obtained must have the form (\ref{series1}), which ,we recall,
is:
\begin{equation}
  \sum_{k \geq 1} \varepsilon^{k}     \left[
       \sum_{\begin{array}{c}
            0 \leq q \leq k \\
            q \equiv k \pmod{2}
             \end{array}}  \cos q\tau   \left(
        \sum_{\begin{array}{c}
              q \leq l \leq k \\
                      l \equiv k \pmod{2}
              \end{array}}   a_{kql} (\lambda /2)^{l} S^{l}(\xi)
                                         \right)
                       \right] \label{series2'}
\end{equation}
where
\[ \xi = \varepsilon x \mbox{ and } \tau = t \sqrt{1-\varepsilon^{2}}.  \]

    This form suggests that there exists a formal solution in
increasing powers of $e^{-\xi}$. The next Theorem shows that such is indeed
the case.

  \begin{theorem}
For every $g$, there exists a formal solution of the corresponding
nonlinear wave equation of the form
\begin{equation}        \label{lseries}
  \sum_{l=1}^{\infty} u_{l}(\tau,\varepsilon) e^{-l\xi},
\end{equation}
with coefficients $u_{l}$ analytic in $\varepsilon$.
The coefficients $u_{l}$ can then be taken rational in $\varepsilon$,
and polynomial in $\cos \tau$.

    The first coefficient $u_{1}$ must have the form
\[      A(\varepsilon) \cos(\tau+\theta)    \]
and all $u_{l}$ can be computed recursively once $A$ and $\theta$ are
given.
  \end{theorem}



    {\em Remarks:}
\begin{enumerate}
    \item  This series has the form (\ref{series2}) announced in the
Introduction. Observe that no conditions on $g$ or parity assumptions
are needed here.
    \item  This series yields, as we shall see, more information
than (\ref{series1}). We have generated its first terms using a
Symbolic Manipulation language (Mathematica). Some interesting points
revealed by this study are developed at the end of the next Section.
    \item  A similar series in increasing powers of $S$ can of
course be constructed, but it is more cumbersome to compute explicitly.
It is again determined by its first term.
    \item  Multiplying the function $A(\varepsilon)$ by an analytic
function which is positive for small real $\varepsilon$ amounts to a
translation in $\xi$.
\end{enumerate}

  \subsection{Proof of Theorem 2.}

    As in the proof of Theorem 1, one sees immediately that each
$u_{l}$ satisfies a differential equation involving the terms previously
computed in the right-hand side. The problem is to show that these
equations are uniquely solvable under our present assumptions.

    Let us write the equation in scaled form:
\[ (1-\varepsilon^{2})u_{\tau\tau}- \varepsilon^{2}u_{\xi\xi}+u
+ g_{2}u^{2}+ g_{3}u^{3}+\cdots = 0.    \]

    Let
\[  M_{l} := (1-\varepsilon^{2})\partial_{\tau}^{2}+(1-l^{2}\varepsilon^{2}). \]

    We then find that $u_{1}$ is determined by the following equation:
\begin{equation}
    (1-\varepsilon^{2})(u_{\tau\tau} + u) = 0.
\end{equation}

    Without loss of generality, we may take
\[ u_{1} = A(\varepsilon) \cos \tau.    \]

    One can see by induction that $A$ will only appear in $u_{l}$
through a factor $A^{l}$. For the sake of simplicity, we therefore
take in what follows $A(\varepsilon) =1$.

    The other coefficients $u_{l}$ satisfy :
\begin{equation}
  M_{l}u_{l} = \sum_{\begin{array}{c}  0 \leq q \leq l  \\
                                       q \equiv l \pmod{2}
             \end{array}}
\alpha_{lq}(\varepsilon) \cos q\tau,
\end{equation}
obtained by setting to zero the coefficient of $e^{-l\xi}$ in (\ref{scaledeq}).

Here,  the right-hand side of the last equation comes from the expansion of
\[ - g_{2}u^{2} - \cdots - g_{l}u^{l}. \]

    It follows that
\begin{equation}
u_{l} = \sum_{l,q}    \alpha_{lq}
         \cos q\tau /(1-q^{2}-\varepsilon^{2}(l^{2}-q^{2})).
\end{equation}

    By induction, one shows that $u_{l}$ is always rational in
$\varepsilon$.

    {\em Remark:} An equation of the form $M_{l}u_{l}=f$ has at most
one solution which is analytic (or even continuous) in $\varepsilon$.
        It is this circumstance which
prevents the occurrence of undetermined terms as in Section 2.
Indeed, such an indeterminacy occurs only for isolated values of
$\varepsilon$, for a given $l$. Note that if
$g$ is odd, there is at most one $2\pi$-periodic solution at all, even without
this assumption.

  \subsection{Relation between the two formal series.}

    As the convergence of the series (\ref{lseries}) will be dealt with in
Section 5, we will now simply point out some connections between the two formal
solutions that we now have at our disposal.

    First of all, as we already remarked, Theorem 1 sets restrictions
on $g$; Theorem 2 does not.

    The existence of the second series clearly amounts to saying
that one can exchange the order of the first two summations in (\ref{series2'}).
But Theorem 2 does not follow from Theorem 1 since it is not clear that
performing the summation  in $k$, which leads to an infinite series,
will leave us with a finite, let alone simply computable, result.

    Theorem 2 does not imply Theorem 1 either, for two different reasons:
\begin{enumerate}
    \item We may always write (\ref{lseries}) as a series in powers
of $S$, but
as before, it is not clear that expanding the $u_{l}$ and
collecting the powers of $\varepsilon$ will leave us with a finite result.
    \item  Even if this were possible, it depends on the choice of
the function $A(\varepsilon)$. To recover (\ref{series1}), we must choose the
expansion of $A$ in such a way that the power of $\varepsilon$ in any term of
the series is not less than the power of $S$ in the same term. One can
always achieve this {\em formally}, but if $A$ thus restricted cannot
be given by a
convergent series, we recover (\ref{series1}), but $A$ is not analytic any
more,  and (\ref{series2'}) doesn't make sense.
\end{enumerate}



\section{Divergence of the First Formal Solution.}

  \subsection{Statement of the result.}

  \begin{theorem}
The series (\ref{series1}) cannot be absolutely convergent on any region of
the form
\[  |\varepsilon| < \frac{1}{\sqrt{2}} + \delta; |e^{-\xi}| < \varrho;
            |\cos \tau| < 1+\varrho     \]
for any choice of the positive numbers $\delta$ and $\varrho$,
and define a breather for these values of $\varepsilon$,
unless $g$ is one of the following:
\[  \sin(\alpha u)/\alpha, \; \sinh(\alpha u)/\alpha, \; u. \]
  \end{theorem}

    {\em Remarks:}
\begin{enumerate}
    \item  The series in question does converge on the domain given in
Theorem 3 for the three exceptional functions we have listed,
as the explicit form
of the corresponding solutions show, provided that $\varrho$ is chosen
small enough with respect to $\delta$. When $g = u$, the series (\ref{series2})
actually reduces to its first term.
    \item Among the three exceptional functions $g$ listed
in the Theorem, only the first
gives rise to breathers. For the other two, the series (\ref{series1})
converges only for $|\xi|$ large enough. Such a behaviour is actually typical,
as is shown by Theorem 4 below.
\end{enumerate}

  \subsection{Proof of Theorem 3.}

    We will first prove that the assumptions of Theorem 3 imply that
the series (\ref{lseries}) converges for large $\xi$ and for {\em complex}
values
of $\varepsilon$ of modulus less than or equal to $\sqrt{2}^{-1}+\delta$.
This, in turn, will lead to the determination of the coefficients of $g$.

    \subsubsection{Reduction to the second series.}

    The hypothesis of the Theorem implies that $u$ is a holomorphic
function of the three variables $\varepsilon$, $e^{-\xi}$, $\cos \tau$ in
the domain indicated. This means that $u$ can also be expanded in
increasing powers of $e^{-\xi}$. We know from Theorem 2 that such a series is
entirely determined by its first term, which, in the notation of
(\ref{series2'}), reads
\[   \sum_{k} a_{k11} \varepsilon^{k} \cos \tau.    \]
Let us call this $a(\varepsilon) \cos \tau$. By assumption, $a$ is analytic,
and is real for real $\varepsilon$. We may take it positive for positive
$\varepsilon$.
It is not identically zero since otherwise $u$ itself would be.

    We claim that $a$ cannot vanish for $\varepsilon \neq 0$.

    Indeed, if $a(\varepsilon_{0}) = 0$, (where $\varepsilon_{0}$ may be
complex),  as the terms of the series
in powers of $e^{-\xi}$ can be found recursively by the procedure of
Section 3, we reach the conclusion that $u(x,t,\varepsilon_{0}) \equiv 0$.
Therefore, $w := \partial u/\partial\varepsilon (x,t,\varepsilon_{0})$
is a solution of:
\[  w_{\tau\tau}+w-\varepsilon_{0}^{2} \Delta w = 0. \]

    Our assumptions on $u$ show that $w$ has a Fourier
expansion $\sum_{j} w_{j} \cos j\tau$, and that for every $j$,
\[  (1-j^{2})w_{j}-\varepsilon_{0}^{2}(\partial_{\xi}^{2}-j^{2})w_{j} =0. \]

    This equation has only one solution tending to $0$ as $|\xi|$ tends
to infinity, namely zero. The argument extends by induction
to all derivatives of $u$ with respect to $\varepsilon$ at
$\varepsilon = \varepsilon_{0}$. As $u$ is assumed to be analytic,
we must have $u \equiv 0$. This contradiction shows that $a$ never vanishes.

    Now $v(\xi - \log a(\varepsilon), \tau, \varepsilon)$ is again
analytic in $e^{-\xi}$ and $\varepsilon \neq 0$, and represents a breather.
As $v$ is a function of $e^{-\xi}$, this solution is single-valued in
$\left\{ 0 < |\varepsilon| < (1/\sqrt{2})+\delta \right\}$.
It has a convergent series of the form (\ref{lseries}), but now with
$u_{1}= \cos \tau$.

    \subsubsection{End of Proof.}

    We now know that (\ref{lseries}) is convergent, with
$u_{1} = \cos \tau$. The coefficients $u_{2}$, $u_{3}$, \ldots can then be
computed recursively. These coefficients
should also be analytic in $\varepsilon$. But from Section 3, we also know that
they are merely rational in $\varepsilon$. $u$ has therefore potential poles
for all the values of
$\varepsilon$ such that there exist $q$ and $l$, nonnegative integers
of the same parity, satisfying
\begin{equation}        \label{poles}
l \geq q \mbox{ and } h(l,q) := 1-q^{2}-\varepsilon^{2}(l^{2}-q^{2}) = 0.
\end{equation}

    We shall express that for these poles not to occur,
some relations must be satisfied by the coefficients of $g$.

    To avoid a pole for $\varepsilon = 1/2$ in $u_{2}$, we must have,
as is found by computing $u_{2}$,
\[  g_{2} = 0.      \]

    There is no condition on $g_{3}$. If it is zero, positive, negative,
we shall see that $g$ must equal $u$, $\sinh(\alpha u)/\alpha$,
$\sin(\alpha u)/\alpha$, with $|g_{3}| = \alpha^{2}/6$.

    Assume $u_{1},\ldots ,u_{l-1}$ have been found, for some $l \geq 4$,
and that they have no pole of modulus $\leq 1/\sqrt{2}$ if and only if
\begin{equation}     \label{excep}
g_{2} = \cdots = g_{[(l-1)/2]} = 0, \; \mbox{and} \;
g_{2m+1} = \frac{(\mbox{sgn}(g_{3}))^{m-1}\alpha^{2m}}
                                  {(2m+1)!}  \end{equation}
for $2m+1 \leq l-1$.

    We claim that this assertion is true with $l$ replaced by $l+1$.



    Let us write the equation giving $u_{l}$ in the form
\begin{equation}    \label{M}
  M_{l}u_{l} = (M_{l}u_{l})_{0}+\varepsilon^{l-1} \cos^{l}\tau
        \,  ((g_{l})_{0}-g_{l}).
\end{equation}

Here, the subscript $0$ refers to the corresponding quantities where
$g$ and $u$ have been replaced by the exceptional values given by
(\ref{excep}).

    Observe now that the last term in (\ref{M})
contributes, if $g_{l} \neq (g_{l})_{0}$, a non-zero term in $\cos 3\tau$
if $l$ is odd, and a constant term if $l$ is even.
Each of them leads to a pole in $u_{l}$, since $(u_{l})_{0}
= M_{l}^{-1}(M_{l}u_{l})_{0}$ has {\em no}
pole. The pole in question occurs for $\varepsilon^{2} = -8/(l^{2}-9)$
for odd $l$, $\varepsilon^{2} = 1/l^{2}$ for even $l$. For $l \geq 4$,
\[ 1/l^{2} < 1/2, \] and for odd $l \geq 4$, \[ 8/(l^{2}-9) \leq 1/2; \]
both types of poles lie in the range $|\varepsilon| <
\sqrt{2}^{-1} + \delta$.

    It follows that $u_{l}$ must have a pole in the forbidden range,
unless $g_{l}=(g_{l})_{0}$,
in which case $u_{l}=(u_{l})_{0}$

    This ends the proof of Theorem 3.

  \subsection{Further remarks.}

    We have seen that the vanishing of the pole corresponding to
$q=0$ or $3$ for every $l$ suffices to determine the coefficients $g_{k}$
one by one. But there are in principle other poles, coming from the higher
harmonics. That all these conditions are automatically satisfied for the
sine-Gordon equation is a remarkable fact. It suggests that the conditions
corresponding to the vanishing of these poles are strongly interrelated.
We list in Table 1 some of these conditions, the more complicated
of which we obtained with the help of a computer. In this Table, we assume
that $g$ is {\em odd}, with
$g_{3} = -1/6$, and we make use of the following abbreviations:
\[  \mbox{(A) : }        1 - 120 g_{5} = 0;              \]
\[  \mbox{(B) : }        1 - 129 g_{5} - 378 g_{7} = 0;  \]
\[  \mbox{(C) : }        g_{5} + 42 g_{7} = 0.           \]

    Observe from the Table that, for instance, the conditions arising
from poles $i$ and $i/\sqrt{5}$ {\em imply} those arising from $i/\sqrt{2}$.
Similarly, those from poles
$i$, $i/3$ and $i/\sqrt{3/2}$ imply that
$g_{5}= 1/120$ or $3/10$;
those from
$i$, $i/3$ and $i/\sqrt{5}$ imply that
$g_{5}= 1/120$ or $589/7680$;
those from
$i$, $i/\sqrt{5}$ and $i/\sqrt{3/2}$ imply that
$g_{5}= 1/120$ or $1/7680$.

    The recurrence of some of the conditions is, however, no surprise:
if a pole has not been cancelled in, say, the computation of $u_{7}$,
it will reappear in $u_{9}$, and the conditions for its disappearance
must contain those omitted in the preceding step of the calculation.

  \begin{table}         \centering
    \begin{tabular}{||l|l|c|r||}        \hline
 $l$ & $q$ & {\em Pole} & \multicolumn{1}{c||}{\em Condition}  \\  \hline
7 & 1 & $i/\sqrt{2} $ & (A)                           \\  \cline{2-4}
  & 3 & $i/\sqrt{2} $ & (A)                           \\  \cline{3-4}
  &   & $i/\sqrt{5} $ & (B)                           \\  \cline{2-4}
  & 5 & $i/\sqrt{2} $ & (A)                           \\  \cline{3-4}
  &   & $i          $ & (C)                           \\  \hline
9 & 3 & $i          $ & (C)                           \\  \cline{3-4}
  &   & $i/\sqrt{2} $ & $(31-2880 g_{5})(1-120 g_{5})=0$ \\  \cline{3-4}
  &   & $i/\sqrt{5} $ & (B)                           \\  \cline{3-4}
  &   & $i/\sqrt{9} $ & $2500-338305 g_{5}+1109760 g_{5}^{2}$ \\
  &   &               & $\mbox{}-1239210 g_{7}-1354752 g_{9} =0 $
                                                       \\  \cline{2-4}
  & 5 & $i          $ & (C)                            \\  \cline{3-4}
  &   & $i/\sqrt{2} $ & $(1-192 g_{5})(1-120 g_{5})=0$ \\  \cline{3-4}
  &   & $i\sqrt{3/7}$ & $125-26985 g_{5}+1739520 g_{5}^{2}$   \\
  &   &               & $\mbox{}+115830 g_{7}+ 746496 g_{9} = 0$
                                                       \\  \cline{3-4}
  &   & $i/\sqrt{5} $ & (B)                            \\  \cline{2-4}
  & 7 & $i\sqrt{3/2}$ & $5 g_{5}+840 g_{5}^{2}+360 g_{7}
                                -10368 g_{9} =0  $     \\  \cline{3-4}
  &   & $i          $ & (C)                            \\  \cline{3-4}
  &   & $i/\sqrt{2} $ & $1-(120 g_{5})^{2} = 0  $      \\  \hline
    \end{tabular}   \caption{Conditions on $g$}
\end{table}

\section{Convergence of the Second Formal Solution.}

  \subsection{Statement of the result.}

    The preceding Section has shown that an expansion of a breather
in increasing powers of $\varepsilon$ cannot converge in the domain in
which the corresponding series for the sine-Gordon equation does. We found
indeed that the solution would then also admit an expansion in powers of
$e^{-\varepsilon x}$, which exhibits poles for finite values of $\varepsilon$.
This fact could not have been found by the study of asymptotics as
$\varepsilon$ tends to zero.

    For the $\varphi^{4}$ model, where
$g(u) = (-(1+u)+(1+u)^{3})/2$, it follows from (19) that poles on the real
axis always occur.

    But poles are all on the imaginary axis if $g$ is
{\em odd}. This suggests that when $g$ is odd, the series (\ref{lseries})
might be convergent, at least for large $\xi$, provided that
$\varepsilon$ is real and {\em kept fixed}. This is precisely the
content of the next Theorem, stated and proved below.

  \begin{theorem}
Assume that $g$ is odd and that its coefficients satisfy:
\begin{equation}    \label{alpha}
       |g_{m}| \leq \frac{\alpha^{m-1}}{m!}
\end{equation}
for all $m$, with some $\alpha > 0$.
Then (\ref{scaledeq}) has a solution given by a series
\[  \sum_{l=1}^{\infty} u_{l}(\tau, \varepsilon)
                              e^{-l\xi} \]
convergent in the same domain as the corresponding series where all
$g_{m}$ are replaced by $\alpha^{m-1}/m!$.
  \end{theorem}

    {\em Remarks:}
\begin{enumerate}
    \item  The majorant series mentioned in the Theorem is
$\alpha u_{SHG}$, where $u_{SHG}$ was defined in \S 1.2, rearranged
after expanding $1/\cosh\xi$ in increasing powers of $e^{-\xi}$.
    \item  It follows from the next Section that there is, up to translation,
only one solution which decays exponentially as $|\xi| \rightarrow \infty$ and
involves odd cosines only. This means that the above function is
{\em the only candidate for an exponentially decaying breather.}
\end{enumerate}

  \subsection{Proof of Theorem 4.}

    We shall prove the Theorem by a majorant method.

\vspace{1 ex}

    {\em Step 1: Preliminaries.}

    We recall that Theorem 3 constructed a formal solution
\[ \sum_{l \geq 1} u_{l} e^{-l\xi}  \]
and that this is equivalent to having a formal solution in powers of $S$,
or $1/\cosh \xi$. We fix now $u_{1} = \varepsilon$ and prove the convergence
of the resulting series.

    If $g$ leads to a breather $u(x,t)$ for some value of the period,
then $\alpha g(u/\alpha)$ admits the breather $\alpha u(x,t)$, and it satisfies
(\ref{alpha}) with $\alpha = 1$. We
therefore assume from now on that $\alpha = 1$.

\vspace{1 ex}

    {\em Step 2: The case $g = \sinh$.}

    We show that by a suitable translation in $\xi$, we can obtain
from $u_{SHG}$ a convergent series solution to
\[ u_{tt} - u_{xx} + \sinh u = 0    \]
with first term $\varepsilon e^{-\xi}\cos \tau$.
    $u_{SHG}$ is equal to
\[ 4i \arctan \left(
\frac{\varepsilon}{\sqrt{1-\varepsilon^{2}}}
\frac{\cos \tau}{i \sinh \xi}     \right)   \]
\begin{eqnarray}
  & = & 4\sum_{l\geq 0} e^{-(2l+1)\xi}
         \sum_{0 \leq s \leq l} \varepsilon^{2s+1}
        \left( \frac{2 \cos \tau}{\sqrt{1-\varepsilon^{2}}}
        \right)^{2s+1}  \nonumber  \\
  & = & 8\sum_{l \geq 0} e^{-(2l+1)\xi}
     \sum_{0 \leq r \leq l} \cos (2r+1)\tau
\nonumber  \\
  &   &  \left[  \sum_{r \leq s \leq l}
        \left( \frac{\varepsilon}{\sqrt{1-\varepsilon^{2}}}
        \right) ^{2s+1}
\left( \begin{array}{c}  2s+1  \\  s-r  \end{array} \right)
\right]  \nonumber  \\
  & := & \sum_{0 \leq s \leq l} b_{ls}(\varepsilon)
                          e^{-(2l+1)\xi} \cos (2s+1)\tau.
\end{eqnarray}

    The series
\[  \sum_{l,s} b_{ls} e^{-(2l+1)\xi} \cos (2s+1)\tau    \]
is absolutely convergent in the domain
\[  \frac{\varepsilon}{\sqrt{1-\varepsilon^{2}}}
    \frac{1}{|\sinh\xi|} < 1.   \]
    As the coefficient of $e^{-\xi}$ in this series is
\[     8\varepsilon\cos \tau \, / \sqrt{1-\varepsilon^{2}}, \]
if we set \[   a(\varepsilon) := \log (\sqrt{1-\varepsilon^{2}}/8), \]
we find
\[ u_{SHG}(\xi-a(\varepsilon),\tau,\varepsilon) = \sum b_{ls}\, e^{-(2l+1)\xi}
(\sqrt{1-\varepsilon^{2}}/8)^{2l+1} \cos (2s+1)\tau \]
and this series now converges for
\[ \left| \frac{\varepsilon/\sqrt{1-\varepsilon^{2}}}
               {\sinh (\xi-a(\varepsilon))}   \right| < 1.  \]
Note that such a region includes the set
\[ \left\{  \xi \geq c, \, 0 \leq \varepsilon \leq d  \right\}  \]
for $1/c$ and $d$ small enough.

    Let now
\[   \bar{a}_{ls}(\varepsilon) := b_{ls}(\sqrt{1-\varepsilon^{2}}/8)^{2l+1}. \]

The series (\ref{lseries}) for $g = \sinh$ is then equal to
\[  \sum_{l,s} \bar{a}_{ls}\, e^{-(2l+1)\xi} \cos (2s+1)\tau.   \]

    The quantities $\bar{a}_{ls}$ are all nonnegative for
$\varepsilon$ nonnegative and fixed.

\vspace{1 ex}

    {\em Step 3: Majorant method and end of proof.}

    Take $g$ as in the Theorem, and consider the series (\ref{lseries})
with, we recall, $u_{1}=\varepsilon$. We shall write it as
\[  \sum_{0\leq s \leq l} a_{ls}(\varepsilon) e^{-(2l+1)\xi}
        \cos (2s+1)\tau.    \]

We have $a_{00} = \bar{a}_{00} =\varepsilon$.
 Let us show by induction on $l$ that
$|a_{ls}| \leq \bar{a}_{ls}$ for $s=0,1,\ldots,l$.

    Assume that the result holds for all integers less than or equal to
$l-1$; $u_{l}$ is then determined by an equation of the form
\[  L_{2l-1}u_{2l-1} = -g_{3} (\ldots)-\cdots
-g_{2l-1} (\ldots). \]

    All terms on the right-hand side are obtained by multiplying
$-g_{q}$, for some $q$, by a linear combination with
{\em nonnegative} coefficients, of terms $a_{ps} \cos (2s+1)\tau$
already computed. To obtain $u_{l}$, one must divide each of them by
$h(2l+1,2s+1)$ as given by (\ref{poles}), since
\[  L_{2l+1} \cos (2s+1)\tau = h(2l+1,2s+1)\cos(2s+1)\tau.  \]

    As $s$ is never greater than $l$, we always have $h < -\varepsilon^{2}$
if $l \geq 3$; in particular, $h$ is negative. Therefore, $|a_{ls}|$ is
no greater than the expression obtained by replacing every $g_{m}$ by
$1/m!$ for odd $m$, and every $a_{ps}$ by $\bar{a}_{ps}$, for $p < l$. The
result of this substitution is precisely $\bar{a}_{ls}$.

    This is the desired result. The formal series solution is
dominated term by term by the corresponding series for the ``sinh-Gordon''
equation, and therefore converges in the same domain.

\section{Decaying solution in $H^{s}$.}

  \subsection{Statement of the result.}

    We assume in this Section that $g$ is an {\em odd} entire function.
We shall prove the following result.

  \begin{theorem}
Let $s$ be a real number greater than $3/2$. For every $T$ greater than, but
close enough to $2\pi$, there exists, up to translation and sign,
a unique $T$-periodic solution of
\[ u_{tt} - u_{xx} + g(u) = 0   \]
which tends exponentially to zero as $x$ tends to $+\infty$, is odd in
$\cos(2\pi t/T)$ and is of class $H^{s}$ in $x$ and $t$.
  \end{theorem}

     The proof will be achieved by the application of the
contraction mapping theorem. The method follows closely the corresponding
argument for ordinary differential equations, with the difference that we
shall follow the dependence of the estimates on $T$ and $s$.

    The next section introduces the notation and defines the map, the
fixed points of which are the desired solutions. The third is devoted to the
action of nonlinear functions on $H^{s}$ spaces. We shall re-derive there a
few classical results; the explicit constants in the inequalities given
here may be new. The two final sections contain the iteration argument and the
proof of our result.

  \subsection{Setting and notations.}

    As in the preceding paragraphs, we introduce the parameter
\[  \varepsilon := \sqrt{1-(2\pi/T)^{2}}    \]
and the scaled variables $\xi$ and  $\tau$. The value of $T$ is fixed for the
rest of this Section.

    We consider the space $X \subset H^{s} \times H^{s-1}$ (for some
$s > 3/2$), which is defined as the set of pairs
\[ \left[ \begin{array}{c}  u \\ v  \end{array} \right] \]
of functions of a variable $\tau$ that are $2\pi$-periodic and have only
odd harmonics in their Fourier expansions. We let
\[ \left\| \left[   \begin{array}{c} u \\ v \end{array} \right] \right\|
_{X}^{2} = |u|_{s}^{2} + \varepsilon^{2}|v|_{s-1}^{2},  \]
\[ A = \left[ \begin{array}{ll}
       0                                                      &  1 \\
\varepsilon^{-2}((1-\varepsilon^{2})\partial_{\tau}^{2} + 1)  &  0
    \end{array}   \right] .     \]
where  we have used the definition
\[  |u|_{s}^2 := \sum_{j} |u_{j}|^{2}(1+j^{2})^{s}      \]
if   $ u = \sum_{j} u_{j} \cos j\tau  $.

    Every element $U$ of $X$ admits the decomposition
\begin{equation}    \label{decomposition}
U = U_{-} + U_{+} + \sum_{j \;\mbox{\scriptsize odd,}\; j \geq 3}
U_{j},
\end{equation}
where
\[  U_{-} = \left[ \begin{array}{c} a \\ -a \end{array} \right] \cos\tau, \;
    U_{+} = \left[ \begin{array}{c} b \\  b \end{array} \right] \cos\tau, \;
    U_{j} = \left[ \begin{array}{c} u_{j} \\ v_{j} \end{array} \right]
\cos j\tau. \]

    One derives easily:
\[  \|U\|_{X}^{2} = \|U_{+}+U_{-}\|_{X}^{2} + \sum_{j} \|U_{j}\|_{X}^{2} \]
and
\begin{equation}  \label{ineq}
\frac{2\varepsilon^{2}}{2+\varepsilon^{2}}
(\|U_{+}\|_{X}^{2} + \|U_{-}\|_{X}^{2})
\leq
\|U_{+}+U_{-}\|_{X}^{2}
\leq
2( \|U_{+}\|_{X}^{2} + \|U_{-}\|_{X}^{2} ).
\end{equation}

    Writing
\[ \omega_{j} = \sqrt{(1-\varepsilon^{2})j^{2} - 1},    \]
we obtain
\[  e^{A\xi} U_{\pm} = e^{\pm\xi} U_{\pm};  \]
\[  e^{A\xi} U_{j}   =
   \left[
    \begin{array}{cc}
   \cos(\omega_{j}\xi/\varepsilon) &
           (\varepsilon/\omega_{j}) \sin(\omega_{j}\xi/\varepsilon) \\
   -(\omega_{j}/\varepsilon) \sin(\omega_{j}\xi/\varepsilon) &
           \cos(\omega_{j}\xi/\varepsilon)
         \end{array}
   \right]  U_{j},
\]
for all odd $j \geq 3$. It follows that
\begin{eqnarray}   \left\|   e^{A\xi}  U_{j}
    \right\| _{X}^{2}  &  \leq  &
C \{ (u_{j})^{2} [(1+j^{2})^{s} +  (1+j^{2})^{s-1}\omega_{j}^{2}]
\nonumber  \\
  &  &  \mbox{} + \varepsilon^{2}(v_{j})^{2}
\left[
(1+j^{2})^{s}/\omega_{j}^{2} + (1+j^{2})^{s-1}
\right]
\} .
\end{eqnarray}



    Now
\[  \left\|
\sum_{j}
U_{j} \right\|_{X}^{2}  =
\sum_j \left( u_{j} \right) ^{2}  (1+j^{2})^{s}  +
\varepsilon^{2} \left( v_{j} \right) ^{2} (1+j^{2})^{s-1} , \]
and if
\[ \varepsilon^{2} \leq \varepsilon_{0}^{2} < 8/9,  \]
there is a positive number $\kappa < 1$ such that
\[  \kappa (1+j^{2}) \leq \omega_{j}^{2} = (1-\varepsilon^{2}) j^{2} - 1
     \leq  1+j^{2}.     \]

    It follows that for every $j$,
\begin{equation}  \left\| e^{A\xi}U_{j} \right\|_{X}^{2}
    \leq (1+1/\kappa) \left\| U_{j} \right\|_{X}^{2}.   \end{equation}

    Observe that this estimate is {\em independent} of $\varepsilon$.

    We also have immediately
\begin{equation}
\left\|
e^{A\xi} U_{\pm}
\right\| _{X}  \leq  e^{\pm\xi} \left\|
U_{\pm}
\right\|  _{X}.
\end{equation}

    We next define an operator $M$ acting on the space
\[ Y \subset C([0,+\infty), X)  \]
consisting of those continuous $X$-valued functions making the quantity
\[ \left\|  \left[
\begin{array}{c}  u  \\  v  \end{array}
   \right]  \right\| _{Y}  : =
\sup_{\xi > 0}    e^{\beta \xi}
\left\|
\left[  \begin{array}{c}  u(\xi)  \\  v(\xi)  \end{array}   \right]
\right\| _{X}   \]
finite; $\beta$ is a given positive number.

    Fix the number $\beta$ between $0$ and $1$, as well as a real number $a$
to be taken sufficiently small later. Define $M$ as follows:
\begin{eqnarray}
M(U)(\xi) &  =  & e^{-\xi}    \left[
\begin{array}{c}  a  \\  -a  \end{array} \right]  \cos \tau    \nonumber \\
          &     & \mbox{} + \int_{0}^{\xi} e^{A(\xi - s)} F_{-}(U)(s) \, ds
\nonumber  \\
          &     & \mbox{} - \int_{\xi}^{\infty} e^{A(\xi - s)}[F_{+}(U)
\nonumber   \\
          &     & \mbox{} + \sum_{j \geq 3} F_{j}(U)](s)\,ds,
\end{eqnarray}
where
\[   F(U) = \left[
\begin{array}{c}  0  \\  g_{3}u^{3}+g_{5}\,\varepsilon^{2}u^{5}
                                +\cdots  \end{array} \right]  =
\left[
\begin{array}{c}  0  \\  h(\varepsilon u) u^{3}  \end{array} \right].  \]
and $F = F_{+} + F_{-} + \sum_{j} F_{j} $ is the decomposition
(\ref{decomposition}) applied to $F(U)$.

{\em Remarks:}

\begin{enumerate}
    \item  The convergence of the improper integral will be shown in
\S 6.4 below; it follows from the fact that $M$ sends $Y$ to itself.
    \item If $M(U) = U$, then writing $U = \left[
\begin{array}{c}  u  \\  v  \end{array} \right] $, one finds that
\[  U_{\xi} = AU + F(U) \]
and therefore that
$\varepsilon u$ solves (\ref{scaledeq}) and that $v = u_{\xi}$. As $u$ has,
by construction, exponential decay, we have achieved the desired solution.
    \item  The uniqueness of the solution will follow from the fact that,
as we prove below, the map $M$ is a contraction near $0$ in $Y$.
\end{enumerate}

    Before we proceed, we shall need a few estimates on products in $H^{s}$.
They are proved in the following section.

  \subsection{Products in $H^{s}$.}

    Let $H^{s}$ denote the ordinary Sobolev space of $2\pi$ periodic
functions of a variable $\tau$. It is well-known that this space is an
algebra for $s > 1/2$. We establish here this fact with an explicit
$s$-dependence of the constants.

    If $u$ has the Fourier expansion
\[  u = \sum_{j= -\infty}^{+\infty} u_{j} e^{ij\tau},   \]
we define its $H^{s}$ norm by
\[  |u|_{s}^{2} = \sum_{j} (1+j^{2})^{s} |u_{j}|^{2}.   \]
(The notation is here slightly different from that in the other
 paragraphs of this paper.)

    We prove our estimates for $C^{\infty}$ functions only, which
as usual, implies the corresponding estimates for non-smooth functions,
by a regularization argument.

    If $u$, $v$ are of class $C^{\infty}$, we have
\[   (uv)_{j} = \sum_{k} u_{k} v_{j-k}.     \]

    As for every $x$, $y$, and positive $\alpha$
\[  (1+y^{2})^{\alpha} + (1+(x-y)^{2})^{\alpha} \geq (1+x^{2}/4)^{\alpha},  \]
we have
\[ (1+j^{2})^{s/2} \leq 2^{s} [ (1+k^{2})^{s/2} + (1+(j-k)^{2})^{s/2}].  \]

    It follows that
\begin{eqnarray}
|(uv)_{j}|^{2}(1+j^{2})^{s} & \leq &
\left(
   \sum_{k} 2^{s} [(1+k^{2})^{s/2} u_{k}v_{j-k} +
       u_{k}(1+(j-k)^{2})^{s/2}v_{j-k}]
\right) ^{2}  \nonumber \\
  & \leq &  2^{2s+1}
\left[
( \sum_{k} (1+k^{2})^{s} u_{k}^{2} ) ( \sum_{k} v_{k}^{2} )
\right.  \nonumber \\
  &  & \left.
\mbox{} + ( \sum_{k} (1+(j-k)^{2})^{s} v_{j-k}^{2}) ( \sum_{k} u_{k}^{2} )
\right] .
\end{eqnarray}

    On the other hand,
\begin{equation}    \label{zeta}
 \sum_{k} u_{k}^{2} \leq |u|^{2}_{s} (\sum_{k} (1+k^{2})^{-s} )
    \leq |u|_{s}^{2}(1+2\zeta(2s)), \end{equation}
with
\[  \zeta (s) = \sum_{n \geq 1} n^{-s}. \]

    If $s > 1/2$, we find
\begin{equation}
|uv|_{s} \leq 2^{s+1} \sqrt{1+2 \zeta (2s) } |u|_{s}|v|_{s}.
\end{equation}

{\em Remark:} One can slightly improve this estimate if $u$ and $v$ contain
only odd harmonics, since the sum in (\ref{zeta}) is taken over odd values
of $k$ only.

\vspace{1 ex}

    We now have, for every $m\geq 1$ and $s$ bounded away from 1/2,
\[      |u^{m}|_{s} \leq (2^{s} C)^{m} |u|_{s}^{m}. \]
It follows that
\begin{equation}
   |g(u)|_{s} \leq \mbox{const.}\sum_{m} |g_{m}|
          (|u|_{s} (2^{s}C))^{m} := \Phi (2^{s}|u|_{s}).
\end{equation}

    If $g(u)$ contains $u^{p}$ as a factor, we obtain
\[  |g(u)|_s \leq |u|_{s}^{p} \Psi_{s} (2^{s/2}|u|_{s}).    \]

  \subsection{Construction of the decaying solution.}

    We show here that  $M(U) = U$ can be solved by the Banach fixed point
theorem in the space $Y$ we have defined.

    To prove this, we estimate, for $U$ and $V$ in $Y$,
\begin{eqnarray}
\| M(U) - M(V) \|_{X} (\xi) & \leq &
      \int_{0}^{\xi} e^{-(\xi -s)} \|F_{-}(U) - F_{-}(V) \|_{X}(s)\,ds
\nonumber  \\
  & &  \mbox{} + \int_{\xi}^{\infty} e^{(\xi - s)}
 \|(F_{+}+\sum_{j} F_{j})(U)   \nonumber \\
  & &   \mbox{} - (F_{+}+\sum_{j} F_{j})(V) \|_{X}(s)\,ds.
\end{eqnarray}

    Now, for every $U$, $\|U\|_{X}$ can be estimated from the
decomposition (\ref{decomposition}) and inequalities (\ref{ineq}).
We thus obtain:
\begin{eqnarray}
\|M(U)-M(V)\|_{Y}  &  \leq  &
C_{1}\varepsilon^{-1}
      \sup_{\xi > 0} \left( \int_{0}^{\xi} e^{-(\xi - s)} e^{-\beta s}
            \|F(U)-F(V)\|_{Y}  \, ds \right.      \nonumber \\
                   &        &
\left.   \mbox{} + \int_{\xi}^{\infty} [1 + e^{(\xi - s)}]  e^{-\beta s}
                        \|F(U)-F(V)\|_{Y} \, ds \right) e^{\beta \xi}
\nonumber \\
                   &   \leq    & C_{1}\varepsilon^{-1}
                                 \|F(U)-F(V)\|_{Y} \sup_{\xi \geq 0}  e^{\beta \xi}
\nonumber  \\
                   &        &      \times
\left(
e^{-\xi} \, \frac{e^{(1-\beta)\xi}}{1-\beta} +
e^{\xi}   \, \frac{e^{-(1+\beta)\xi}}{1+\beta}
\right) ,
\end{eqnarray}
so that
\[  \|M(U)-M(V)\|_{Y}  \leq
\frac{2C_{1}/\varepsilon}{1-\beta^{2}}\|F(U)-F(V)\|_{Y} . \]

    Now, by Section 6.2, if we write
\[  U = \left[ \begin{array}{c} u \\ u' \end{array} \right] , \]
\[  V = \left[ \begin{array}{c} v \\ v' \end{array} \right] , \]
we have
\begin{eqnarray}
\|F(U)-F(V)\|_{X} & = & \varepsilon  |u^{3} h(\varepsilon u) -
v^{3} h(\varepsilon v) |_{s-1}  \nonumber \\
                  & = &  \varepsilon
\left| \int_{0}^{1} f(u+t(v-u))(v-u) \, dt \right|_{s-1} ,
\end{eqnarray}
where
\[ f(\xi) := d(\xi^{3} h(\varepsilon \xi))/d\xi . \]

    It follows that, using the results of \S 6.3,
\begin{eqnarray}
\|F(U)-F(V)\|_{X} & \leq & \varepsilon  C_{2} (|u|_{s-1},|v|_{s-1})
                                    |v-u|_{s-1} 2^{s}  \nonumber  \\
                  & \leq & \varepsilon  C_{3}(s,\|U\|_{X},\|V\|_{X})
                \|U-V\|_{X},
\end{eqnarray}
where $C_{3}$ is an increasing function of its arguments.
       As $ \| \;\;\; \|_{X} \leq \| \;\;\; \|_{Y} $,
\begin{equation}
\|M(U)-M(V)\|_{Y} \leq (2/(1-\beta^{2}))
C_{1}C_{3}(s,\|U\|_{Y},\|V\|_{Y})   \|U-V\|_{Y} .
\end{equation}

    Observe now that the map $M$ is a contraction on a sufficiently
small ball in $Y$, uniformly in $\varepsilon$,
because $f(0) = 0$ ensures that $C_{3}$ tends to zero with
$\|U\|_{Y} + \|V\|_{Y}$. If therefore
\[   M(0) = \left[
\begin{array}{c}  a  \\  -a  \end{array}
\right] \cos \tau       \]
is small enough, the iterates
$(M^{n}(0))_{n\geq 1}$ converge as $n \rightarrow \infty$ to a
fixed point of $M$.

    This ends the construction of the decaying solution.



    \subsection{Uniqueness.}

    We prove here the uniqueness part of the Theorem. From the results
of the previous section, we know that there is one solution in a suitably
small neighborhood of the origin in $Y$. This solution clearly does not
depend on the choice of the number $\beta$ introduced in \S 6.4, but it might
depend on $a$. We show here that different values of $a$ correspond to
the translates of the solution obtained in \S 6.4.

    Let us write $M_{a}$ for $M$, to stress its dependence on $a$.

    If $M_{a}(U) = U$, where $U$ lies in $Y$ for some value of $\beta$,
then the translates $U_{\lambda}=U(.+\lambda)$ tend to zero in $X$ (and $Y$),
and therefore $V_{\lambda} = U_{\lambda}-M_{0}(U_{\lambda})$ also tends
to zero. But as $U_{\lambda}$ and $M_{0}(U_{\lambda})$ both solve
\[    dU/d\xi = A U + F(U_{\lambda}),   \]
we have
\[    dV_{\lambda}/d\xi = A V_{\lambda} \]
using the definition of $M_{a}$. As $V_{\lambda}$ decays exponentially,
it follows that
\[   V_{\lambda} = \left[
\begin{array}{c} a(\lambda) \\ -a(\lambda) \end{array}
     \right] \cos \tau .    \]
Furthermore, $a(\lambda)$ is continuous and
tends to zero as $\lambda$ tends to infinity.
It follows that for $\lambda$ large enough, $V_{\lambda}$
is the unique fixed point of $M_{a(\lambda)}$. As $U \not\equiv 0$,
$a(\lambda)\neq 0$ for large $\lambda$, and $a(\lambda)=a$ has a solution
for $a$ small enough and of the same sign as $a(\lambda)$.

    This proves the uniqueness, up to sign and translation, of the small,
odd, exponentially decaying solutions of the problem.

\section{Appendix: Bibliographical Remarks.}

    We here briefly discuss the literature on the problem. We have limited
ourselves to the works that have a direct relation to our results.
 The contributions are grouped according to the method they use.

    1) The authors of \cite{DHN,kosevich-kovalev} and others first observed that
 the first terms of a formal expansion in powers of $\varepsilon$
could be found for the $\varphi^{4}$-model, in which one takes
\[  g(u) = \frac{1}{2}(-(1+u) + (1+u)^{3}).  \]
This suggested the existence of breathers for
an equation other than sine-Gordon.
    This result is explained by our Theorem 1. Our second formal
expansion, however, is not to be found in the literature.
 Endeavors have been made to combine this result
with a Fourier expansion with a view to obtain a contradiction at the formal
level: In \cite{elonsky,elonsky-silin}, by truncating the Fourier series and using
numerical integration of the remaining equations, one tries to extract
information by ``matching'' a solution tending to zero at $+\infty$
with one tending to zero at $-\infty$. It is found there that the coefficients of
$g$ should be restricted, but that the sine does not satisfy
these restrictions. In \cite{kruskal-segur}, an attempt is made to use the method of
matched asymptotic expansions to produce a correction to the solution
of order $e^{-C/|\varepsilon|}/\varepsilon$, such as
would mean that the solution
cannot be even in $x$, despite the formal solution being even
to all orders. Numerical integration of approximate solutions is also used,
together with several truncations
where in particular $2\pi$ is substituted for the period $T$ ---  a case
in which (see below) breathers such as considered there do not exist.

    Both of these arguments are motivated by phenomena that are
found in other problems. The first is reminiscent of calculations of the
splitting of separatrix loops in some ordinary differential equations,
where in particular the successes of Melnikov's method are well-known.
The second is motivated by the relevance of ``transcendentally
small'' corrections to the solutions of some recent problems.

    2) A number of attempts at finding breathers numerically
have been made (see e.g.
\cite{ablowitz-kruskal-ladik,eilbeck,makhankov,kudryavasev,wingate}). One way to
generate such solutions is to investigate a ``head-on'' collision of two
``kink'' solutions (solutions having different limits as $x$ goes to
plus or minus infinity). Such kinks are easy  to write down for the
$\varphi^{4}$-model, since time-independent kinks exist. One finds that
way long-lived ``quasi-breather'' modes. They appear to fade away after
some time.

    3) Energy estimates can be used to obtain bounds on possible
breathers, as well as a limited number of nonexistence results.
 We mention here \cite{coron} which proves under mild decay assumptions
 that breathers, if they exist,
must have period greater than or equal to $2\pi$, \cite{vuillermot2} which
shows nonexistence in a number of cases (all these cases have,
in our notation, $\lambda \leq 0$, so that even formal solutions do not
exist in these cases). McKean (unpublished) also obtained some uniqueness
results by related methods; he also investigated the possibility of
finding a ``separated'' solution $f(h_{1}(x)h_{2}(t))$, as in the
sine-Gordon case; he finds that no other wave equation has breathers of
that form.

    4) It is natural to try to perturb the sine-Gordon breather
itself \cite{weinstein-birnir}. If one sets
\[    g(u) = \sin u +\alpha h(u),   \]
and tries a perturbation expansion in powers of $\alpha$,
\[  u = u_{SG} + \alpha u_{1} + \cdots, \]
one finds that $u_{1}$ cannot be found unless $h$ satisfies some
conditions. Concrete (nonexistence) results in that direction are in
preparation (\cite{weinstein-birnir}).

    5) Attempts have been made to view the problem as a dynamical
system with $x$ as ``time variable'', whereby the problem becomes one of
finding an orbit homoclinic to zero. We have mentioned formal arguments in that
direction earlier. A. Weinstein observed that the stable manifold theorem
gives, for every value of the period greater than, and close to, $2\pi$,
a three-dimensional manifold of solutions tending to zero as $x$ goes to
$+\infty$ (but not necessarily $-\infty$). The dimension of this stable
manifold increases with the period. It contains possible  breathers, together
with time-independent solutions. Several authors, see e.\ g. \cite{vuillermot1},
using classical
ideas from dynamical systems, show that solutions dominated by
$c/\cosh x$, with $c$ small enough must vanish identically. Such results
apply to the sine-Gordon equation as well: the constant $c$ {\em depends}
on the period.




\end{document}